\documentclass[11pt]{article}

\usepackage{amsmath,amssymb}

\usepackage{graphics}
\usepackage{tikz}

\newcommand{\beq}[1]{\begin{equation}\label{#1}}
\newcommand{\enq}{\end{equation}}
\newcommand{\raf}[1]{(\ref{#1})}
\newcommand{\qed}{\hfill\hbox{\rule{6pt}{6pt}} \medskip}
\newcommand{\proof}{\noindent {\bf Proof}.\ }

\newcommand{\cB}{\mathcal{B}}

\def\gal{\alpha}
\def\gbe{\beta}

\def\ZZ{\mathbb{Z}}

\DeclareMathOperator*{\argmax}{argmax}

\DeclareMathOperator*{\argmin}{argmin}

\newtheorem{theorem}{Theorem}
\newtheorem{lemma}[theorem]{Lemma}

\newtheorem{corollary}[theorem]{Corollary}

\newtheorem{definition}[theorem]{Definition}

\newtheorem{remark}[theorem]{Remark}

\def\problemDefInstance{Instance}

\def\problemDefQuestion{Question}

\def\problemDef#1#2#3{
\begingroup
\offinterlineskip
\centerline{\vbox{
\halign{\vrule \hskip1em ## \hfil & ## \hfil \hskip1em\vrule\cr
\noalign{\hrule}
\multispan2\vrule height1.3em depth.8em \hfil{\sc #1}\hfil\vrule\cr
\noalign{\hrule}
\multispan2\vrule height.8em \hfil \vrule \cr
{\bf \problemDefInstance{} :} &
\bgroup
\multiply\hsize by 2
\divide\hsize by 3
\parindent0pt
\baselineskip1.2em
\lineskip1pt
\lineskiplimit0pt
\vtop{#2}
\egroup\cr
\multispan2\vrule height.8em \hfil \vrule \cr
{\bf \problemDefQuestion{} :} &
\bgroup
\multiply\hsize by 2
\divide\hsize by 3
\parindent0pt
\baselineskip1.2em
\lineskip1pt
\lineskiplimit0pt
\vtop{#3}
\egroup\cr
\multispan2\vrule height.4em \hfil \vrule \cr
\noalign{\hrule}
}
}}
\endgroup
\vskip.5em
} 

\usepackage{mathrsfs}

\title{Recognizing distributed approval voting forms and correspondences.}
\author{
Endre Boros\thanks{
MSIS Dep. of RBS and RUTCOR, Rutgers University,
100 Rockafeller Road, Piscataway, NJ 08854-8054, USA.
(endre.boros@rutgers.edu)}
\and
Ond\v{r}ej \v{C}epek\thanks{
Department of Theoretical Informatics and Mathematical Logic, Charles University,
Malostransk\'{e} n\'{a}m. 25, 11800 Praha 1, Czech Republic.
(ondrej.cepek@mff.cuni.cz)}
\and
Vladimir Gurvich\thanks{
National Research University Higher School of Economics, Moscow, Russia.
( vgurvich@hse.ru , vladimir.gurvich@gmail.com)}
\and
Kazuhisa Makino\thanks{
Research Institute for Mathematical Sciences $($RIMS$)$ Kyoto University,
Kyoto 606-8502, Japan. (makino@kurims.kyoto-u.ac.jp)}
}
\date{\today{}}

\begin{document}
\maketitle

\begin{abstract}
We consider distributed approval voting schemes.
Each voter  $i \in I$  has  $\alpha_i$  cards that
(s)he  distributes among the candidates  $a \in A$ as a  measure of approval.
One (or several) candidate(s) who received the maximum number of cards is
(are) elected.
We provide polynomial algorithms
to recognize voting forms and  voting correspondences generated by such
voting schemes in cases when either the number of candidates or the number of voters is equal to $2$.
We prove that for two voters, if $\alpha_2\geq \alpha_1-2\geq 0$ 
then the unique voting correspondence has distinct rows. We also characterize voting forms with distinct rows.
\end{abstract}

\section{Introduction}
Let $I = \{1, \ldots, n\}$ be a set of $n$ voters (electors, players) and    
$A = \{a_1, \ldots, a_p\}$  be a set of $p$ candidates
(decisions, outcomes). In the literature, several voting schemes are considered in which
one (or all) candidate(s) who received the maximum number of votes is (are) elected.

In {\em plurality voting}
(see for example, \cite{BCG19,DL04,Mei15,MPRJ17})  each voter chooses only one candidate.
In {\em approval voting} each voter can choose (approve) an arbitrary subset  $A' \subseteq A$  of candidates.
This model was considered by Brams and Fishburn \cite{BF78}; see also \cite{Alo06,BF07}.
Masso and M. Vorsatz \cite{MV06} suggested further generalization, called \emph{weighted approval voting}, waving the neutrality condition:
each voter has a positive integer weight (which may correspond, for example, to the number of people that (s)he
represents). Then, naturally, the elected candidate(s) is (are) the candidate(s) with the maximum  weighted sum of approval.

\smallskip 

In this paper we focus on \emph{distributed approval voting}: Each voter $i \in I$  is given a positive integer number $\alpha_i$
of voting cards that (s)he can distribute arbitrarily among the candidates, thus expressing the measure of approval by this voter.
This notion generalizes both approval and weighted approval voting. 
The candidates with the highest number of voting cards are chosen for the final round, and eventually one of them is selected. 
The integer $\alpha_i$ is called the {\em weight} or {\em power} of voter $i$ and
the obtained voting scheme is called {\em distributed approval voting}. In this paper we consider both the sets of candidates entering the final round (distributed approval voting correspondences) and the possible final round selections (distributed approval voting forms). 

The standard mathematical approach in voting theory (see, for example, \cite{Mou83,Pel84}) is
to consider the normal form of a voting scheme. In this paper we also use normal form to represent distributed approval voting. 
Most of our results in this paper
are related to the case of $n=2$ voters. We therefore adopt a simplified notation where the number of voting
cards of the two voters is denoted by  $\gal$ and $\gbe$  instead of $\alpha_1$ and $\alpha_2$. 
In Section \ref{p=2} we focus on the case of two candidates and arbitrary number of voters, 
so there we switch back to the notation $\alpha_i$, $i\in I$ for the powers of voters. 

For $n=2$ we denote by
\[
X ~=~ \left\{ x\in \ZZ_+^A ~\left|~ \sum_{a\in A} x_a=\gal \right.\right\}
~~~\mbox{and}~~~ Y ~=~ \left\{ y\in \ZZ_+^A ~\left|~ \sum_{a\in A} y_a=\gbe \right.\right\}
\]
the sets of possible voting strategies of the two voters. We denote by
\begin{equation}\label{e0}
k=\binom{\gal+p-1}{p-1}=|X| ~~~\mbox{ and }~~~ \ell=\binom{\gbe+p-1}{p-1}=|Y|
\end{equation}
the cardinalities of these sets.

To simplify our notation, for voting strategies $z\in \ZZ_+^A$ we also introduce 
\[
AM(z)=\argmax_{a\in A} z_a
\]
to denote the subset of candidates that receive the most voting cards under voting strategy $z$. 

A matrix $h\in \left(2^A\right)^{k\times \ell}$ is called a \emph{voting correspondence}.
We call it a \emph{distributed approval voting correspondence} 
(\emph{DAV correspondence} for short) 
if we can label the sets of voting strategies $X=\{x^1, \dots ,x^k\}$ and $Y=\{y^1, \dots ,y^\ell\}$ such that
\[
h(i,j) ~=~ AM(x^i+y^j)
\]
for all $i=1, \dots ,k$ and $j=1, \dots ,\ell$.
We call a matrix $g\in A^{k\times \ell}$ a \emph{voting form}.
We call it a \emph{distributed approval voting form} (\emph{DAV form} for short) 
if we can label the sets of voting strategies $X=\{x^1, \dots ,x^k\}$ and $Y=\{y^1, \dots  ,y^\ell\}$ such that
\[
g(i,j) ~\in~ AM(x^i+y^j)
\]
for all $i=1, \dots ,k$ and $j=1, \dots ,\ell$.

As a small example assume that $A=\{a,b\}$ and $\gal=\gbe=3$.
Then we have $X=Y=\{(3,0),(2,1),(1,2),(0,3)\}$, and $k=\ell=4$. 
See Figures \ref{fig1} and \ref{fig2} for the corresponding DAV correspondence and some possible DAV forms. 

\begin{figure}[htb]
\[
\begin{array}{||c||c|c|c|c||}
\hline\hline
&(3,0)&(2,1)&(1,2)&(0,3)\\
\hline\hline
(3,0)&\{a\}&\{a\}&\{a\}&\{a,b\}\\
\hline
(2,1)&\{a\}&\{a\}&\{a,b\}&\{b\}\\
\hline
(1,2)&\{a\}&\{a,b\}&\{b\}&\{b\}\\
\hline
(0,3)&\{a,b\}&\{b\}&\{b\}&\{b\}\\
\hline
\hline
\end{array}
\]
\caption{A DAV correspondence with $p=2$ and $\gal=\gbe=3$.\label{fig1}}
\end{figure}

\begin{figure}[htb]

\[
\begin{array}{||c||c|c|c|c||}
\hline\hline
&(3,0)&(2,1)&(1,2)&(0,3)\\
\hline\hline
(3,0)&a&a&a&a\\
\hline
(2,1)&a&a&b&b\\
\hline
(1,2)&a&b&b&b\\
\hline
(0,3)&b&b&b&b\\
\hline
\hline
\end{array}
~~~~~~~~~
\begin{array}{||c||c|c|c|c||}
\hline\hline
&(3,0)&(2,1)&(1,2)&(0,3)\\
\hline\hline
(3,0)&a&a&a&b\\
\hline
(2,1)&a&a&b&b\\
\hline
(1,2)&a&a&b&b\\
\hline
(0,3)&a&b&b&b\\
\hline
\hline
\end{array}
\]
\caption{Two DAV forms with $p=2$ and $\gal=\gbe=3$. The fist one has no identical row, while the second one has identical rows. \label{fig2}}
\end{figure}

Note that while a given set of parameters $p$, $\gal$, and $\gbe$
defines a unique DAV correspondence,
we may have many different DAV forms corresponding to it.
In some of these DAV forms we may have identical rows,
see Figures \ref{fig1} and \ref{fig2}.

Recall that for given $p$, $\alpha$, and $\beta$ we can easily compute $k$ and $\ell$ by \eqref{e0}. On the other hand, given a voting correspondence or voting form, we can determine $p$ just by inspecting the matrix entries (since all candidates must appear in any DAV correspondence or form). Consequently we can also determine the values of $\alpha$ and $\beta$ or conclude that no $\alpha$ and/or $\beta$ can satisfy \eqref{e0}. In the rest of the paper we assume that for the considered voting correspondences or forms the parameters $p$, $\alpha$, and $\beta$ that satisfy \eqref{e0} always exist. We also assume that if a voting correspondence or form is given, then the set $A$ of candidates is simply the set of elements appearing in the given matrix. 

The paper is structured as follows. Section 2 deals with DAV correspondences for two voters ($n=2$). It contains two main results: sufficient and necessary conditions for any DAV correspondence satisfying these conditions to have distinct rows (or columns, by symmetry), and a polynomial time recognition algorithm which for a given input voting correspondence decides whether it is a DAV correspondence, and in the affirmative case outputs a corresponding labelling of rows and columns by voting strategies. Section 3 deals with DAV forms for two voters. Similarly, it also contains two main results. The first one provides sufficient and necessary conditions for any DAV form satisfying these conditions to have distinct rows (or columns, by symmetry). Here the set of conditions is more complex than for DAV correspondences, and also the proof is more involved. The second result is a polynomial time recognition algorithm which for a given input voting form decides whether it is a DAV form, and in the affirmative case outputs a corresponding labeling of rows and columns by voting strategies. However, this algorithm assumes additional conditions on the input parameters $p$, $\alpha$ and $\beta$, and therefore, unlike for correspondences, it cannot be used for every input voting form. Section 4 treats three special cases not covered completely by the previous sections. The first special case $\alpha = \beta = 1$ is the original plurality voting (for two voters) which motivated the concept of distributed approval voting. For this case we give a complete characterization of DAV voting forms by a set of three forbidden submatrices. For the second special case $\alpha, \beta \leq 2$ we give a simple polynomial time recognition algorithm for DAV forms as this case is not fully covered by the general recognition algorithm of Section 3. The last special case deals with two candidates ($p=2$) and an arbitrary number of voters ($n\geq 2$). We conclude the paper by several remarks and open problems in Section 5.

\section{Distributed Approval Voting Correspondences}

Our first result characterizes the cases when all
corresponding DAV forms have pairwise distinct rows.

\begin{theorem}
\label{th-2a}
For $p\geq 2$ and $\gal\geq 2$ a DAV correspondence has distinct rows if and only if $\gbe\geq \gal-2$.
\end{theorem}

\proof
Let us show first that $\gbe<\gal-2$ implies that the DAV correspondence has two identical rows. Consider two voting strategies $x,x' \in X$ with $x_a = \gal$ (and hence $x_b=0$ for all $b \neq a$) and $x'_a = \gal-1$. Since $\gbe<\gal-2$ candidate $a$ received enough voting cards in both $x$ and $x'$ so that $AM (x+y)=AM (x'+y) = \{a\}$ for all $y\in Y$. Hence the DAV correspondence has identical rows $x$ and $x'$.

Let us assume next that $\gbe\geq \gal-2$ and for two arbitrary distinct voting strategies $x,x'\in X$ consider the following two cases:

\begin{description}
\item[\rm Case 1:] $AM(x) \neq AM(x')$. Without loss of generality we assume $a\in AM(x)\setminus AM(x')$. Choose $b\in AM(x')$. This choice of $a$ and $b$ implies that $x_a \geq x_b$ and $x'_a < x'_b$.
Let us consider voting strategy $y \in Y$ which splits all $\gbe$ voting cards among $a$ and $b$ (i.e. $y_c=0$ for all $c\in A\setminus\{a,b\}$) in such a way that $y_b\leq y_a\leq y_b+1$. If $\gbe$ is even, then $y_a = y_b$ and hence $x_a + y_a \geq x_b + y_b$ while $x'_a + y_a < x'_b + y_b$. Therefore we have $a\in AM(x+y)\setminus AM(x'+y)$. If $\gbe$ is odd, then $y_a = y_b +1$ and hence $x_a + y_a > x_b + y_b$ while $x'_a + y_a \leq x'_b + y_b$. Therefore we have $b\in AM(x'+y)\setminus AM(x+y)$. In either case we proved that $AM(x+y)\neq AM(x'+y)$.

\item[\rm Case 2:] $AM(x)=AM(x')$. Choose $a\in AM(x)$ and without loss of generality assume that $x_a\geq x'_a$. If $x_b\geq x'_b$ for all $b \in A$,  then $x\geq x'$, which in turn implies $x = x'$ (because the sum of all coordinates is $\gal$ for both $x$ and $x'$) which is a contradiction. Hence there exists a candidate $c \in A$ such that $0 \leq x_c < x'_c$ and therefore the set $B=\{c\in A\mid x'_c-x_c > x'_a-x_a\}$ is not empty. Let us now fix an arbitrary $b\in \argmax_{c\in B} x'_c$. Clearly $x'_b >0$. Note also that $x_b<x_a$ since $B\cap AM(x)=\emptyset$ by the definition of $B$ and the assumption $a \in AM(x)=AM(x')$.

Let us show that there exists a voting strategy $y \in Y$ such that all $\gbe$ voting cards are split among candidates $a$ and $b$ in such a way that $x'_a+y_a\leq x'_b+y_b\leq x'_a+y_a+1$ holds. This is not difficult to verify since $\gbe\geq \gal-2$ is assumed and $x'_a > x'_b > 0$ implies $x'_a < \gal$. If $x'_a+y_a=x'_b+y_b$ then we have $b\in AM(x'+y)\setminus AM(x+y)$, while if $x'_a+y_a+1=x'_b+y_b$ then we have $a\in AM(x+y)\setminus AM(x'+y)$. In either case we have $AM(x+y)\neq AM(x'+y)$.
\end{description}
Thus the rows corresponding to voting strategies $x$ and $x'$ are different. Since this is true for two arbitrary rows, the statement follows.
\qed

Now we proceed to the main result of this section which is a polynomial time recognition algorithm for DAV correspondences. We start with introducing additional notation.

For a voting strategy $x\in X$ and a candidate $a\in A$ let us denote by $s(x,a)$ the number of occurrences of $a$ in the row of the DAV correspondence labeled by $x$, i.e.,

\[
s(x,a) ~=~ |\{y\in Y\mid a\in AM(x+y)\}|.
\]
We call the vector $s(x)=(s(x,a)\mid a\in A)$ the \emph{signature} of $x$.

\begin{lemma}\label{l0}
If $\gbe\geq\gal-1$, then the rows of a DAV correspondence have pairwise distinct signatures.
\end{lemma}

\proof
Let us consider arbitrary two distinct voting strategies $x,x'\in X$, and set $B= \argmax_{a\in A}(x_a-x'_a)$.
We first claim that for every $y\in Y$ the containment $AM(x'+y) \cap B \subseteq AM(x+y)$ holds.
Indeed,  for any $a \in A$ and $b \in AM(x'+y)$, we have $x'_b+y_b \geq x'_a+y_a$.
Moreover, if $b \in AM(x'+y) \cap B $ then $x_b+y_b \geq x_a+y_a$ and thus $b \in AM(x+y)$.

The above containment implies that $s(x',b) \leq s(x,b)$ for every $b \in B$. To complete the proof we show
that for at least one $b \in B$ this inequality is strict, i.e., $s(x',b) < s(x,b)$.
\begin{description}
\item[\rm Case 1:]  $p=2$.  In this case,  $|B| = |A\setminus B| = 1$. Let $B = \{b\}$ and $A\setminus B = \{a\}$. By definition $x_b > x'_b$ and $x_a < x'_a$. Since $\gbe\geq\gal-1$ there exists $y \in Y$ for which $x_a+y_a\leq x_b+y_b\leq x_a+y_a+1$. From these inequalities it is not hard to derive that $b \in AM(x+y) \setminus AM(x'+y)$, impliying that $s(x',b) < s(x,b)$.

\item[\rm Case 2:]   $p\geq 3$.   Define $b\in \argmax_{d\in B} x'_d$ and $c\in \argmax_{d\in A\setminus B} x'_d$. Note that $\max(x'_b,x'_c)$ is equal to the maximum of the components of $x'$.
Note also that we have $x'_c>0$ since otherwise $x\geq x'$ follows, contradicting that $x$ and $x'$ are distinct elements in $X$.
Since $\gbe\geq \gal-1$ and $p \geq 3$, there exists a $y\in Y$ such that $x'_b+y_b+1=x'_c+y_c$ and $y_b+y_c\geq \gbe-1$.
Then we have $b\in  AM(x+y) \setminus AM(x'+y)$ by $b \in B$.
Thus  $s(x',a)<s(x,a)$.\qed
\end{description}

Given a voting correspondence $h\in \left(2^A\right)^{k\times \ell}$ we define $s_h(i,a)$  as the number of candidates $a\in A$ that occur in row $i$ of $h$:
\[
s_h(i,a) = |\{j\in \{1,\dots,\ell\}\mid a\in h(i,j)\}|.
\]
We call the vector $s_h(i)=(s_h(i,a)\mid a\in A)$ the \emph{signature of row} $i$ of $h$, or simply a \emph{row-signature} of $h$.

\begin{theorem}\label{t-corr-rec}
Let $h\in \left(2^A\right)^{k\times \ell}$ be a voting correspondence.
Then we can recognize in polynomial time if $h$ is a DAV correspondence.
If so, we can also obtain in polynomial time  a  labeling of  rows and columns by the elements of $X$ and $Y$, respectively. Moreover,  if $\ell\geq k$ then the labeling of the rows is unique, otherwise the labeling of the columns is unique.
\end{theorem}

\proof
Let us recall that if $h$ is indeed a DAV correspondence then we can compute the unique $\alpha$ and $\beta$ values. Without loss of generality we can assume that $\ell \geq k$ since otherwise we can interchange the roles of $k$ and $\ell$ and transpose the input matrix. Thus, $\beta\geq \alpha$ can also be assumed. 
Let us next compute the row-signatures of $h$ and compare them to the signatures of the voting strategies $x\in X$.
If $s_h(i)=s(x)$, then we assign $x$ to the $i$th row of $h$, i.e., we label $x=x^i$.
By Lemma \ref{l0},  $h$ is a DAV correspondence if and only if we have exactly one such $x \in X$ for every row of $h$. 
If it is the case, we obtain in $O(k\ell p)$ time a unique assignment of $X=\{x^1,\dots,x^k\}$ to the rows of $h$.

Next we check if $Y$ can be assigned to columns of $h$. For this we create a bipartite graph $G=(V,E)$ with $V=Y\cup\{1,\dots, \ell\}$ and  $E=\{(y,j) \mid AM(x^i+y)=h(i,j)$ for all $i=1, \dots ,k\}$ in $O(\ell^2kp)$ time.
It is not difficult to see that $G$ has a perfect matching if and only if  $h$ is a DAV correspondence,  where the matching provides a one-to-one assignment of $Y$ to the columns of $h$.  Since a perfect matching can be computed in $O(\ell^3)$ time if exists \cite{Egervary31,Konig31,Kuhn55},
the labeling of columns can be computed in $O(\ell^3p)$ time. Since such a perfect matching may not be unique, our labeling of the columns may not be unique either. 
\qed

\section{Distributed Approval Voting Forms}

In this section, we consider DAV forms.
For voting strategies $x,x'\in X$ let us define
\[
D(x,x') ~=~ \{y\in Y\mid AM(x+y)\cap AM(x'+y)=\emptyset\}
\]
as the set of voting strategies $y\in Y$ that surely differentiate the rows labeled by voting strategies $x$ and $x'$ of any DAV form.

\begin{lemma}\label{l1}
Given distinct voting strategies $x,x'\in X$ choose $a\in AM(x-x')$ and $b\in AM(x'-x)$ and consider the voting strategy $x''\in X$ defined by
\[
x''_e ~=~ \begin{cases}
x'_a+1 & e=a,\\
x'_b-1 & e=b,\\
x'_e & e\in A\setminus\{a,b\}.
\end{cases}
\]
Then we have
\[
D(x,x'') \subseteq D(x,x').
\]
\end{lemma}

\proof
For  $x,x'\in X$ with $x\neq x'$, let $a\in AM(x-x')$ and $b\in AM(x'-x)$.
We assume that there exists a $y\in D(x,x'')\setminus D(x,x')$, and derive a contradiction.

By definition,  we have $AM(x+y)\cap AM(x''+y)=\emptyset$ and $AM(x+y)\cap AM(x'+y)\neq \emptyset$.
Hence, we have  a candidate $c\in (AM(x+y)\cap AM(x'+y))\setminus AM(x''+y)$, and the following inequalities are satisfied:
\begin{subequations}
\begin{align}
x_c+y_c &\geq x_e+y_e ~~~\mbox{  for all } e\in A,\label{1a}\\
x'_c+y_c&\geq x'_e+y_e ~~~\mbox{  for all } e\in A,\label{1b}\\
x''_c+y_c &< x''_d+y_d ~~~\mbox{  for some } d\in AM(x''+y).\label{1c}
\end{align}
\end{subequations}
By definition, $a\not=b$ and $c\not=d$ hold.
Moreover, \eqref{1b} and \eqref{1c} imply that $c=b$ or $d=a$.
In the following we separately consider the following three cases.
\begin{description}
\item[\rm Case 1.]
$c=b$ and $d=a$: Recall that by the choice of $a$ and $b$ we have
\begin{equation}
\label{eq-ae11}
x_a-x'_a>0>x_b-x'_b.
\end{equation}
Therefore,  we have
\[
\begin{split}
x_a+y_a ~=~ (x_a-x'_a-1) +x''_a+y_a  &>~ (x_a-x'_a-1) +x''_b+y_b ~=~ (x_a-x'_a-2) + x'_b+y_b\\ &\geq~ (x_b-x'_b)+x'_b+y_b ~=~ x_b+y_b,
\end{split}
\]
where the first and second inequalities follow from \eqref{1c} and \eqref{eq-ae11}, respectively.
This contradicts \eqref{1a}.

\item[\rm Case 2.] $c=b$ and $d\in A\setminus\{a,b\}$: In this case $x'_d=x''_d$, and thus \eqref{1b} and \eqref{1c} imply $x'_b+y_b=x'_d+y_d$. Since $x_b-x'_b\leq x_d-x'_d$ by our choice of $b$,   we also get $x_d+y_d\geq x_b+y_b$, which together with \eqref{1a} and $c=b$ implies that $d\in AM(x+y)$. This contradicts $AM(x+y)\cap AM(x''+y)=\emptyset$.

\item[\rm Case 3.]
$c\in A\setminus\{a,b\}$ and $d=a$: In this case we have $x'_c=x''_c$ and thus inequalities \eqref{1b} and \eqref{1c} imply that $a\in AM(x'+y)\cap AM(x''+y)$. By our choice of $a$ we also have $x_a-x'_a\geq x_e-x'_e$ for all $e \in A$. Therefore, $a\in AM(x'+y)$  implies $a\in AM(x+y)$, which contradicts $AM(x+y)\cap AM(x''+y)=\emptyset$.
\end{description}
\qed

Lemma~\ref{l1} implies that decreasing the distance $\sum_{i=1}^p |x_i - x'_i|$ between two (row) voting strategies $x$ and $x'$ (i.e. moving from $x'$ to $x''$) can only shrink the set of (column) voting strategies that surely differentiate the two rows. Hence, in order to verify that every pair of (row) voting strategies $x$ and $x'$ has nonempty $D(x,x')$, it suffices to verify this condition for "neighboring" pairs, i.e. those  $x$ and $x'$ for which $x_a=x'_a+1$, $x_b=x'_b-1$, and $x_e=x'_e$ if $e\not\in \{a,b\}$ for some candidates $a$ and $b$.

\begin{lemma}\label{l-ade2}
For $p \geq 3$, let $x$ and $x'$ be two voting strategies in $X$ such that $x_a=x'_a+1$, $x_b=x'_b-1$,
and $x_e=x'_e$ if $e\not\in \{a,b\}$ for some candidates $a$ and $b$. Let $y \in Y$ be arbitrary, Then
$AM(x+y) \cap AM(x'+y)=\emptyset$ implies $AM(x+y)=\{a\}$ and  $AM(x'+y)= \{b\}$.
\end{lemma}

\proof
We first note that $b \not\in AM(x+y)$, since otherwise also $b \in AM(x'+y)$ holds, which contradicts the disjointness assumption $AM(x+y) \cap AM(x'+y)=\emptyset$. Similarly we have $a \not\in AM(x'+y)$.

We next assume by contradiction that there exists a candidate $c \in AM(x+y)$ such that $c\not\in \{a,b\}$.
Let us select $d  \in AM(x'+y)$. If $d\not\in \{a,b\}$, then we must have $x_d=x_c=x'_d=x'_c$, which implies $c \in AM(x'+y)$, a contradiction to the disjointness assumption. By $c \in AM(x+y)$, we have $x'_c=x_c \geq x_a > x'_a$ and so $d\not=a$ follows. Hence $AM(x'+y)=\{b\}$ is the only remaining possibility.  However, in this case $x'_b > x'_c$ and so $x_b = x'_b-1 \geq x'_c = x_c$ implying $b \in AM(x+y)$, which again contradicts the disjointness assumption. Therefore there exists no $c\not\in \{a,b\}$ such that $c \in AM(x+y)$, and by a symmetrical argument also no $c\not\in \{a,b\}$ such that $c \in AM(x'+y)$.

Combining the arguments in the above two paragraphs proves the claim, i.e. $AM(x+y)=\{a\}$ and  $AM(x'+y)= \{b\}$.
\qed

\begin{theorem}
\label{th-eerty1}
Given $p$, $\gal$ and $\gbe$, all distributed approval
voting forms have distinct rows if and only if
one of the following conditions hold:
\begin{description}
\item[\rm (i) $p=2$:]  $\gbe\geq \gal-1$ and $\gal+\gbe$ is odd;
\item[\rm (ii) $p=3$:]   $\gbe\geq 2\gal$ or $\gbe=2\gal-2$;
\item[\rm (iii) $p \geq 4$:]  $(\gal,\gbe)\not=(1,1)$ and $\gbe\geq 2\gal -2$.
\end{description}
\end{theorem}

\proof
We separately prove each of the three cases (i), (ii), and (iii).

(i). For $p=2$ let us denote $A=\{a,b\}$.
Note that if $\gal+\gbe$ is odd, we always have a unique elected candidate for every pair of voting strategies $x \in X$ and $y \in Y$.
Therefore, in this case the DAV correspondence represents in fact a unique DAV form. Moreover, note that $\gbe = \gal-2$ is impossible for
$\gal+\gbe$ odd, and so $\gbe\geq \gal-2$ is equivalent to $\gbe\geq \gal-1$. Thus, by Theorem \ref{th-2a}, if  $\gal+\gbe$ is odd,
the unique DAV form has distinct rows if and only if $\gbe\geq \gal-1$.

Let us now consider the case in which  $\gal+\gbe$ is even. We show in this case that there always exists a DAV form
with two identical rows.  Assume first that both $\gal$ and $\gbe$ are even. Then we consider two voting strategies $x,x' \in X$ where
$x=(\gal/2,\gal/2)$ and $x'=(\gal/2-1,\gal/2+1)$. Let $y=(\gbe_1, \gbe_2) \in Y$ be arbitrary.
If $\gbe_1\geq \gbe_2+2$,  we have $a \in AM(x+y) \cap AM(x'+y)$.
If $\gbe_1\leq \gbe_2$, we have $b \in AM(x+y) \cap AM(x'+y)$.
Note that $\gbe_1\not= \gbe_2+1$ by the assumption that $\gbe$ is even.
We next assume that both $\gal$ and $\gbe$ are odd. Then we consider voting strategies $x=(\lceil\gal/2\rceil,\lfloor\gal/2\rfloor)$, $x'=(\lfloor\gal/2\rfloor,\lceil\gal/2\rceil)$, and $y=(\gbe_1, \gbe_2)$. If $\gbe_1\geq \gbe_2+1$,  we have $a \in AM(x+y) \cap AM(x'+y)$.
If $\gbe_1\leq \gbe_2-1$, we have $b \in AM(x+y) \cap AM(x'+y)$.
Note that  $\gbe_1\not= \gbe_2$  by the assumption that $\gbe$ is odd.
\medskip

(ii). For $p=3$ let us denote $A=\{a,b,c\}$.
Let us first show that if  $\gbe\geq 2\gal$ or $\gbe=2\gal-2$ than any DAV form has pairwise distinct rows .
Let $x$ and $x'$ be  two voting strategies in $X$ such that $x_a=x'_a+1$, $x_b=x'_b-1$ and $x_c=x'_c$.
We prove that  the rows corresponding to $x$ and $x'$ are distinct which by Lemma \ref{l1} implies that any two rows are distinct.
To this end, we construct a voting strategy $y = (y_a,y_b,y_c)\in Y$ for which we prove that $AM(x+y)=\{a\}$ and $AM(x'+y)=\{b\}$.
Let us define $(y_a,y_b,y_c)$ such that
\begin{equation}
\label{eq-def}
x'_a+y_a=x_b+y_b ~{\rm and} ~ y_c \in \{0,1\}.
\end{equation}
Voting strategy $y\in Y$ satisfying (\ref{eq-def}) can be constructed as follows. We distribute the $\gbe$ voting cards in such a way that we first add $|x'_a - x_b| \leq \gal -1 \leq \gbe$ cards to either $y_a$ or $y_b$ depending which of $x'_a, x_b$ is smaller, and then add pairs of cards to $y_a$ and $y_b$ until all $\gbe$ cards are used except maybe one, which is then given to $y_c$. Note that since up to $\gbe \geq 2(\gal-1)$ cards were distributed into $y_a$ and $y_b$, the construction of $y$ guarantees $x'_a+y_a \geq \gal -1$ and $x_b+y_b \geq \gal -1$ which implies
\begin{equation}
\label{eq-rect}
x_a+y_a \geq \gal ~{\rm and} ~ x'_b+y_b \geq \gal.
\end{equation}
Furthermore note that
\begin{equation}
\label{eq-ew1}
x_c=x'_c\leq \gal - 1, ~ x_a+y_a > x_b+y_b,  ~{\rm and} ~ x'_a+y_a < x'_b+y_b,
\end{equation}
where the first inequality follows from $x_a \geq 1$ and the remaining two are an easy consequence of (\ref{eq-def}).

Let us consider the case of $\gbe=2\gal-2$. If $x_c=x'_c= \gal - 1$, then $x=(1,0,\gal - 1)$, $x'=(0,1,\gal - 1)$, and $y=(\gal - 1,\gal - 1,0)$.
Hence we have $AM(x+y)=\{a\}$ and $AM(x'+y)=\{b\}$.
On the other hand, if $x_c=x'_c \leq \gal - 2$, then we have  $x_a+y_a > x_b+y_b$ by (\ref{eq-ew1})
and  $x_a+y_a \geq \gal > x_c+y_c$ by (\ref{eq-rect}) implying $AM(x+y)=\{a\}$.
Similarly, $x'_b+y_b > x'_a+y_a$ holds by (\ref{eq-ew1})
and  $x'_b+y_b \geq \gal > x_c+y_c$ by (\ref{eq-rect}) implying $AM(x+y)=\{b\}$.

Let us consider the case of $\gbe \geq 2\gal$. In this case we can replace (\ref{eq-rect}) by stronger inequalities $x'_a+y_a \geq \gal$
and $x_b+y_b \geq \gal$ which imply $x_a+y_a\geq \gal +1 > x_c+y_c$ and
$x'_b+y_b \geq \gal +1 > x'_c+y_c$. These together with \raf{eq-ew1} again imply that  $AM(x+y)=\{a\}$ and $AM(x'+y)=\{b\}$.

We next consider the case in which  $\gbe= 2\gal-1$ or $\gbe \leq 2\gal-3$. Let $x$ and $x'$ be two voting strategies such that $x=(1,0,\gal-1)$ and $x'=(0,1,\gal-1)$.
We show that $AM(x+y) \cap AM(x'+y)\not=\emptyset$ holds for any $y \in Y$, which completes the proof of (ii).
By contradiction assume that there exists $y \in Y$ such that $AM(x+y) \cap AM(x'+y)=\emptyset$.
By Lemma \ref{l-ade2}, we have $AM(x+y) =\{a\}$ and $AM(x'+y)=\{b\}$.
Thus $y$ satisfies $y=(\lambda,\lambda, \mu)$ for some nonnegative integers $\lambda$ and $\mu$ and moreover we may assume $\mu\in \{0,1\}$.
If  $\gbe= 2\gal-1$, we have $\lambda=\gal-1$ and $\mu=1$ which implies that $c\in AM(x+y) \cap AM(x'+y)$.
If $\gbe \leq 2\gal-3$, we have $\lambda\leq \gal-2$ which again implies $c\in AM(x+y) \cap AM(x'+y)$.
Thus, in either case, we derive a contradiction.
\medskip

(iii). For $p\geq 4$ let us denote $A=\{a,b,c,d \ldots\}$. This case is somewhat similar to (ii) so we make the proof less detailed.
Let us first show that if  $(\gal,\gbe)\not=(1,1)$ and $\gbe\geq 2\gal -2 $ then any DAV form has pairwise distinct rows.
Let $x$ and $x'$ be  two voting strategies in $X$ such that $x_a=x'_a+1$, $x_b=x'_b-1$ and $x_e=x'_e$ for all $e\not\in \{a,b\}$.
We prove that  the rows corresponding to $x$ and $x'$ are distinct which by Lemma \ref{l1} implies that any two rows are distinct.
To this end, we construct a voting strategy $y \in Y$ for which we prove that $AM(x+y)=\{a\}$ and $AM(x'+y)=\{b\}$.
Let $c$ be a candidate in $A$ such that $c \in \argmin_{e \in A \setminus \{a ,b\}}x_e$ and let us define $y$ such that
$x'_a+y_a=x_b+y_b$, $y_c \in \{0,1\}$, and $y_e=0$ if $e\not\in \{a,b,c\}$. We can again derive both (\ref{eq-rect}) and (\ref{eq-ew1}) with the only difference that this time $x_c=x'_c\leq \frac{\gal - 1}{2}$ because there are at least two other candidates outside of $\{a,b\}$ and $c$ is the one with fewer (or equal number of) voting cards. Furthermore note that
\begin{equation}
\label{eq-ew2}
x_a+y_a \geq \gal > x_e+y_e  ~{\rm for~all}~ e \in A\setminus \{a,b,c\}, ~{\rm and} ~ x'_b+y_b \geq \gal > x'_e+y_e  ~{\rm for~all}~ e \in A\setminus \{a,b,c\}.
\end{equation}
Now if $\gal>1$, we have
\[
x_a+y_a \geq \gal > \frac{\gal-1}{2} + 1 \geq x_c+y_c  ~{\rm and}~  x'_b+y_b \geq \gal > \frac{\gal-1}{2} + 1 \geq x'_c+y_c,
\]
which together with (\ref{eq-ew1}) and (\ref{eq-ew2}) implies $AM(x+y)=\{a\}$  and $AM(x'+y)=\{b\}$.
If $\gal=1$ and $\gbe\geq 2$, then  we have $x=(1,0,0 \ldots)$, $x'=(0,1,0 \ldots)$, and $y=(\lambda,\lambda \ldots)$ for $\lambda\geq 1$, and so
\[
x_a+y_a \geq 2 > 1 \geq x_c+y_c  ~{\rm and}~  x'_b+y_b \geq 2 > 1 \geq x'_c+y_c,
\]
which again together with (\ref{eq-ew1}) and (\ref{eq-ew2}) implies $AM(x+y)=\{a\}$  and $AM(x'+y)=\{b\}$. This proves our claim.

We next  consider the case in which  $(\gal,\gbe)=(1,1)$ or $\gbe \leq 2\gal-3$. Let $a,b,c \in A$,  and
let $x$ and $x'$ be two strategies such that $x_a=1$, $x'_a=0$, $x_b=0$, $x'_b=1$, $x_c=x'_c=\gal-1$,
and $x_e=x'_e=0$ if $e\not\in \{a,b,c\}$.
We show that $AM(x+y) \cap AM(x'+y)\not=\emptyset$ holds for any $y \in Y$, which completes the proof of (iii).
By contradiction assume that there exists $y \in Y$ such that $AM(x+y) \cap AM(x'+y)=\emptyset$.
By Lemma \ref{l-ade2}, we have $AM(x+y) =\{a\}$ and $AM(x'+y)=\{b\}$.
Thus $y$ satisfies $y_a=y_b=\lambda$.
If $(\gal,\gbe)=(1,1)$, then $\lambda=0$. If $d$ is the candidate with $y_d=1$ then  $d \in AM(x+y) \cap AM(x'+y)$.
If $\gbe \leq 2\gal-3$ then $\lambda\leq \gal-2$ which implies $c\in AM(x+y) \cap AM(x'+y)$.
Therefore, in either case, we derive a contradiction.
\qed

In the rest of this section we study the problem of recognizing if a given voting form $g\in A^{k\times \ell}$ is distributed approval. Recall that it is easy to check equalities \eqref{e0}, and thus we can assume those in the sequel. To formulate our result, we need a few more definitions.

Given a subset $B\subseteq A$, we denote by
\[
s_g(i,B) ~=~ |\{j\mid g(i,j)\in B, j=1,\dots,\ell\}|
\]
the number of candidates belonging to $B$ in row $i$ of $g$.
For a voting strategy $x\in X$ we denote by
\[
L(x,B) ~=~ \{ y\in Y\mid AM(x+y)\subseteq B\} ~~~\mbox{and}~~~ U(x,B) ~=~ \{ y\in Y\mid AM(x+y)\cap B\neq \emptyset\}
\]
the sets of voting strategies of the second (column) voter for which we must (resp. can) have a candidate from $B$ in the row of a DAV form labeled by $x$. Note that $L(x,B)\subseteq U(x,B)$ always hold. Note also that $x$ can label row $i$ of voting form $g$ only if
\begin{equation}\label{e-L-U}
|L(x,B)| \leq s_g(i,B)  \leq |U(x,B)|
\end{equation}
holds for all subsets $B\subseteq A$.

Let us consider two voting strategies $x,x'\in X$ and associate to them the set
\[
B(x,x') ~=~ \{a\in A\mid x_a>x'_a\}.
\]
By definition, for any distinct $x$ and $x'$ in $X$ we have $B(x,x')\neq \emptyset$ and $B(x,x') \cap  B(x',x)=\emptyset$.

\begin{lemma}\label{l-main}
If $p\geq 3$ and $\gbe\geq 2\gal$, then for all $x,x'\in X$, $x\neq x'$ we have $|U(x',B(x,x'))|<|L(x,B(x,x'))|$.
\end{lemma}
\proof
Let us denote $B=B(x,x')$.
Note that $x_b>x'_b$ for all $b\in B$, and $x'_c\geq x_c$ for all $c\in A\setminus B$.
Hence if $y\in Y$ satisfies that $AM(x'+y) \cap B\not=\emptyset$, then we have $AM(x+y) \subseteq B$,
implying that  $U(x',B)\subseteq L(x,B)$. Thus, to prove the claim, it remains to show that $L(x,B)\setminus U(x',B) \not=\emptyset$.

Choose two candidates $b$ and $c$ such that $b\in \argmax_{a\in B} x_a$ and $c\in \argmax_{a\in B(x',x)} x'_a$ and  construct $y\in Y$ such that $x_b+y_b=x'_c+y_c$ and $y_b+y_c\geq \gbe -1$. 
If $\gbe$ and $x_b-x_c'$ are of different parity then we must have $y_a=1$ for some $a\not\in\{b,c\}$. Such an $a$ exists because we assume $p\geq 3$. 
Since $x_b\geq 1$, $x'_c\geq 1$,  and $\gbe\geq 2\gal$,  we get that $x_b+y_b=x'_c+y_c\geq \gal+1$, while $y_a+\max\{x_a,x'_a\}\leq \gal$ for every $a\not\in \{b,c\}$.
It is now not difficult to verify that $AM(x+y)=\{b\}$ and $AM(x'+y)=\{c\}$ which implies $y\in L(x,B)\setminus U(x',B)$.
\qed

Note that for $p=3$ and $\gbe=2\gal -1$ the above claim is not true, e.g., for voting strategies $x=(1,0,\gal-1)$ and $x'=(0,1,\gal-1)$.
A similar claim can be shown for $p=2$ if $\gbe\geq \gal-1$ and $\gal + \gbe$ is odd, which is equivalent to the condition in Theorem \ref{th-eerty1} (i).

Now we are ready to prove our main result for the recognition of DAV forms.

\begin{theorem}\label{t2}
Let $g$ be a voting form such that $p\geq 3$ and $\gbe\geq 2\gal$ hold for the corresponding parameters. 
Then we can recognize in polynomial time if $g$ is DAV form. 
If so, we can also obtain in polynomial time a labeling of rows and columns by the elements of $X$ and $Y$, respectively. 
Moreover, the labeling of rows is unique.
\end{theorem}

\proof
As before, we assume that equalities \eqref{e0} hold, otherwise $g$ is not a DAV form.
Let us next define
\[
\cB ~=~ \{B(x,x')\mid x,x'\in X\}.
\]
Then for each row $i$ of $g$ let us check the inequalities \eqref{e-L-U} for all $B\in\cB$ and $x\in X$. By Lemma \ref{l-main} there can be at most one $x\in X$ such that \eqref{e-L-U} hold for all $B\in\cB$. If to all rows there is one such $x\in X$ corresponding, then $g$ is distributed approval,
and we constructed a unique assignment of the voting strategies in $X$ to the rows of $g$.

If such a unique row labeling exists, then we can follow the proof of Theorem \ref{t-corr-rec} and construct a bipartite graph between the elements of $Y$ and the columns of $g$. We define $G=(V,E)$ with $V=Y\cup\{1,\dots, \ell\}$ and  $E=\{(y,j) \mid AM(x^i+y)\ni g(i,j)$ for all $i=1, \dots ,k\}$. It is again easy to see that any perfect matching in $G$ provides us with a labeling of the columns of $g$ by the elements of $Y$. 
\qed

Note that if $p\geq 3$ and $\alpha\geq 2\beta$, then we can repeat the previous claim for the columns instead of the rows. 

\section{Special Cases}

In this section we shall investigate DAV forms in three special cases, namely $\alpha = \beta =1$, $\alpha, \beta \leq 2$ (in both cases with $n=2$ as before), and $p=2$ with arbitrary $n$. For $\alpha = \beta =1$ we give a characterization of DAV forms by a collection of forbidden submatrices. For $\alpha, \beta \leq 2$ (which includes the previous case as a subcase) we give a recognition algorithm for DAV forms which uses similar techniques as Theorem~\ref{t2} (note that the cases $\alpha = \beta =1$ and $\alpha = \beta =2$ are not covered by Theorem~\ref{t2}). Finally, for $p=2$ (only two candidates) we derive a characterization of DAV forms for any number of voters $n$. 

\subsection{Case $\alpha = \beta =1$.}




This case corresponds to what is known as plurality voting with two voters\footnote{A generalization, called a separable discrete function, was considered in \cite{BCG19}.}. 
Here we get $k=\ell=p$ by \eqref{e0}, so every DAV form is a square $p \times p$ matrix where each of the $p$ candidates serves once as a row label (voting strategy) and once as a column label because the sets $X$ and $Y$ of row and column voting strategies are in this case just sets of all unit vectors of dimension $p$ that can be  identified with the set $A$ of candidates.  Hence recognizing whether a given voting form $g$ is a DAV form amounts to finding the two appropriate permutations of $A$ (one for rows and one for columns) or proving that such permutations do not exist. 
Let us define the following three matrices
\begin{equation}
\label{forbidden}
m_1 = \left[ {\begin{array}{cc}
a & b \\
c & a
\end{array} } \right]
\hspace{2cm}
m_2 = \left[ {\begin{array}{cc}
a & a \\
a & a
\end{array} } \right]
\hspace{2cm}
m_3 = \left[ {\begin{array}{cccc}
a & a & b & b
\end{array} } \right]
\end{equation}
where candidate $a$ is distinct from both candidates $b$ and $c$ (but $b=c$ is allowed).

\begin{lemma}
\label{two-pairs}
Let $g$ be a $p \times p$ voting form, $i \neq j$ two distinct parallel lines (two rows or two columns) in $g$, and let $a$ be a candidate such that both $i$ and $j$ contain at least two occurrences of $a$ each. Then $g$ contains $m_1$ or $m_2$ as a submatrix.
\end{lemma}

\proof
Assume that lines $i,j$ are rows (the proof for columns is completely analogous). Take one occurrence of $a$ in row $i$, say in column $k$ (that is $g(i,k)=a$), and one occurrence of $a$ in row $j$ in some column $\ell \neq k$ ($g(j,\ell)=a$). Such entries in $g$ of course exist by the assumption of the lemma. Now if $g(i,\ell)=g(j,k)=a$ then $g$ contains $m_2$. So assume without loss of generality that $g(j,k) \neq a$. Now if also $g(i,\ell) \neq a$ we have that $g$ contains $m_1$. So let us assume $g(i,\ell)=a$. By assumption of the lemma row $j$ contains two occurrences of $a$, so let $g(j,m)=a$ for some column $m$ different from both $k$ and $\ell$. Now if $g(i,m)=a$ we have $m_2$ in rows $i,j$ and columns $\ell,m$, while if $g(i,m) \neq a$ we have $m_1$ in rows $i,j$ and columns $k,m$, which finishes the proof.
\qed

\begin{lemma}
\label{too-few}
Let $g$ be an $p \times p$ voting form which contains at most $p-1$ candidates. Then $g$ contains $m_1$ or $m_2$ as a submatrix.
\end{lemma}

\proof
By the cardinality constraint, each row $i$ of $g$ contains some entry $a_i$ twice (there may be more than one such entry, in such a case pick any of them as $a_i$). Since there are $p$ rows and only $p-1$ distinct $a_i$'s, there must be two distinct rows $i \neq j$ such that $a_i = a_j$. Such two rows fulfill the assumptions of Lemma~\ref{two-pairs} and so the claim follows.
\qed

\begin{definition}
\label{greedy}
Let $g$ be an $p \times p$ voting form containing candidates from set $A$ of size $p$. Let us define a greedy assignment (GA for short), which assigns to each row and column of $g$ a (possibly empty) subset of $A$ in the following manner: $a \in A$ is assigned to line $j$ if and only if line $j$ contains at least two occurrences of $a$.
\end{definition}

\begin{lemma}
\label{partial-permutation}
Let $g$ be an $p \times p$ voting form containing candidates from set $A$ of size $p$, such that $g$ contains no $m_1$ or $m_2$ or $m_3$ as a submatrix. Then GA assigns at most one label to every row and column, and both the assignment to rows and the assignment to columns are injective, i.e. no two rows are assigned the same label and no two columns are assigned the same label.
\end{lemma}

\proof
Since $g$ does not contain $m_3$ as a submatrix, no line (row or column) in $g$ gets assigned more than one label. Now assume by contradiction, that the same label $a$ was assigned to two distinct rows $i \neq j$. That means by the definition of the greedy assignment that both $i$ and $j$ contain a pair of $a$'s, and hence rows $i$ and $j$ fulfill the assumptions of Lemma~\ref{two-pairs}. This gives a contradiction with the assumption, that $g$ contains neither $m_1$ nor $m_2$. The argument for columns is identical.
\qed

We have proved so far that for a voting form $g$ with no $m_1$ or $m_2$ or $m_3$ as a submatrix procedure GA in fact assigns "partial" permutations of $A$ to both rows and columns of $g$. Now we shall prove that these "partial" permutations fulfill the conditions required by a DAV form and can be easily extended to full permutations.  

\begin{lemma}
\label{can-be-extended}
Let $g$ be an $p \times p$ voting form containing candidates from set $A$ of size $p$, such that $g$ contains no $m_1$ or $m_2$ or $m_3$ as a submatrix. Then the labeling created by GA can be extended to permutations of $A$ for both rows and columns which define a DAV form.
\end{lemma}

\proof
There are two cases to consider.
\begin{enumerate}
\item All entries in $g$ are covered by the labels assigned by GA, i.e. for all $1 \leq i,j \leq p$ we have that $g(i,j) = a$ implies that row $i$ was assigned label $a$ by GA, or column $j$ was assigned label $a$ by GA (or both). In this case we take the injective assignment for rows given by GA (injectivity is guaranteed by Lemma~\ref{partial-permutation}) and complete it in an arbitrary manner to a permutation of all $p$ labels, and we do the same thing for columns.

\item There exists an entry $x$ in row $i$ and column $j$ of $g$ covered neither by the label assigned by GA to row $i$ nor by the label assigned by GA to column $j$. It is obvious from the definition of GA, that there is no other occurrence of $x$ in row $i$ or column $j$, which in turn implies that there is no other occurrence of $x$ in the entire voting form $g$, as such an occurrence, say in row $k$ and column $\ell$, would give $m_1$ in the intersections of rows $i$ and $k$ with columns $j$ and $\ell$. Thus $x$ has a single occurrence in $g$.

Now let us assume by contradiction there there is another entry $y \neq x$ in $g$ which is also not covered by GA. By the same arguments as above, we can conclude that also candidate $y$ has a single occurrence in $g$. Now deleting row $i$ and the column containing $y$ we arrive to an $(n-1) \times (n-1)$ voting form which contains at most $(n-2)$ distinct candidates (all occurrences of both $x$ and $y$ were deleted). Now by Lemma~\ref{too-few} we get that $g$ must contain $m_1$ or $m_2$ as a submatrix, which is a contradiction. Thus all other entries in $g$ except of the single occurrence of $x$ are covered by GA.

The fact that $x$ occurs only once in $g$ implies that every row other than $i$ and every column other than $j$ contain at most $n-1$ distinct candidates, which means that at least one candidate is repeated in every line except of $i$ and $j$, and hence GA assigns a label to every line except of $i$ and $j$ (and of course no such label may be $x$). Therefore, by Lemma~\ref{partial-permutation}, all rows except of $i$ are assigned a permutation of the $n-1$ remaining labels except of $x$, and the same is true for columns except of $j$. Since no label can be assigned to two rows or two columns by Lemma~\ref{partial-permutation}, the only possibility is that both row $i$ and column $j$ are not assigned any label by GA. So to complete the labels given by GA to full permutations, it suffices to assign candidate $x$ to both $i$ and $j$.
\qed
\end{enumerate}

Lemma~\ref{can-be-extended} immediately implies the following corollary.

\begin{corollary}
\label{sufficient}
Let $g$ be an $p \times p$ voting form containing candidates from set $A$ of size $p$, such that $g$ contains no $m_1$ or $m_2$ or $m_3$ as a submatrix.  Then $g$ is a DAV form.
\end{corollary}

Let us now observe, that the sufficient condition used in Corollary~\ref{sufficient} is also necessary.

\begin{lemma}
\label{necessary}
Let $g$ be an $p \times p$ voting form containing candidates from set $A$ of size $p$, such that $g$ contains one of $m_1$ or $m_2$ or $m_3$ as a submatrix.  Then $g$ is not a DAV form.
\end{lemma}

\proof
If $g$ contains $m_1$ then the two $a$ entries in it cannot be covered by a single line and hence any assignment covering both of them must use one row and one column with label $a$ (it cannot use two distinct rows or two distinct columns with the same label). But now the entry at the intersection of this row and column which is distinct form $a$ is not covered. Thus $g$ is not a DAV form. If $g$ contains $m_2$ then it is immediately obvious the the four $a$ entries in it cannot be covered by one row and one column with label $a$. Thus $g$ is not a DAV form. Finally, if $g$ contains $m_3$ then again it is quite obvious that if the line of $g$ containing $m_3$ is assigned label $a$, then there is no way to cover the two $b$ entries (only one column may be assigned label $b$), and vice versa. Thus it follows again that $g$ is not a DAV form.
\qed

We are now ready to state the main result of this section. 
\begin{theorem}
\label{characterization}
Let $g$ be an $p \times p$ voting form containing candidates from set $A$ of size $p$. Then $g$ is a DAV form if and only if $g$ contains no $m_1$ or $m_2$ or $m_3$ as a submatrix. Moreover, $\{m_1, m_2, m_3\}$ is a unique minimal set with this property.
\end{theorem}
\proof
Corollary~\ref{sufficient} and Lemma~\ref{necessary} proves the first half of our claim. 

What remains to show, is that the set $\{m_1, m_2, m_3\}$ is minimal, that is, if we leave any one of the three matrices out, then there exists $p \times p$ voting form $g$ with exactly $p$ distinct candidates, which is not a DAV form, and does not contain any of the remaining two $m_i$'s as a submatrix. Let us consider the following three voting forms
\begin{equation}
g_1 = \left[ {\begin{array}{ccc}
a & b & b \\
c & a & b \\
c & c & a
\end{array} } \right]
\hspace{2cm}
g_2 = \left[ {\begin{array}{ccc}
a & a & a \\
a & a & b \\
a & a & c
\end{array} } \right]
\hspace{2cm}
g_3 = \left[ {\begin{array}{cccc}
a & a & b & b \\
a & c & c & b \\
a & c & d & b \\
a & d & d & b
\end{array} } \right]
\end{equation}
which are all of type $p \times p$ with exactly $p$ distinct candidates. Clearly each $g_i$ contains the corresponding $m_i$ and so none of the $g_i$'s is a DAV form. Moreover, no $g_i$ contains any of the other two $m_j$'s, $j \neq i$, as a submatrix, proving the minimality of the set $\{m_1, m_2, m_3\}$.
\qed

Note that recognizing if a given voting form is DAV can be done by simply running procedure GA in $O(p^2)$ time. 
According to the above results, if GA does not produce a valid labeling then the given form is not DAV, 
while if it does, then it is DAV.

\begin{remark}
The game forms that appear in plurality voting for two voters can be further generalized. Such generalizations were recently considered in several papers
\cite{BCG19,BGMP10,MPRJ17,Kuk11} under the name of {\em separable} or {\em assignable} game forms. 
In this generalization the dimensions of the voting form can be arbitrary and the same candidate may appear as a label of rows and/or columns arbitrary number of times. 
Some interesting relations between improvement acyclicity, Nash-solvability, and assignability were found.
\end{remark}

\subsection{Case $\alpha \leq 2$, $\beta \leq 2$, and $p \geq 3$.}
In this case we do not give a full characterization by a set of forbidden submatrices as in the $\alpha = \beta = 1$ case, instead we provide a polynomial time recognition algorithm for all voting forms of this type. If $\alpha = \beta = 1$ then  the polynomial time recognition algorithm is a byproduct of the characterization theorem (Theorem~\ref{characterization}). If $\alpha = 1$, $\beta = 2$ or $\alpha = 2$, $\beta = 1$ then the recognition algorithm follows directly from Theorem~\ref{t2}. Hence the remaining case is $\alpha = 2$, $\beta = 2$ (and $p \geq 3$) which we shall solve using counting arguments, where for each row we count the number of particular entries in that row. This is very similar to the concept of row signatures used in the lemmas proving Theorem~\ref{t2}, in particular to the use of inequalities \eqref{e-L-U} where we will consider only one-element subsets. We proceed by case analysis for all possible types of row and column labels, each such label being a pair of candidates (not necessarily distinct). 

The row labeled by $aa$ has at least $p(p-1)/2 + 1 = (p^2 - p + 2)/2$ occurrences of $a$ (at all columns labeled by distinct pairs and at column labeled by $aa$) and at most $p(p+1)/2$ occurrences of $a$ (all columns may contain $a$ in this case). For a fixed $b \neq a$ the row labeled by $aa$ has no forced occurrence of $b$  and at most $1$ occurrence of $b$ (the possible occurrence being at the column labeled by $bb$).

The row labeled by $ab$ for $a \neq b$ has at least $p-1$ occurrences of $a$ ( at $p-1$ columns labeled by pairs of entries containing $a$ and not containing $b$) and at most $p-1 + (p-2)(p-3)/2 + 1 = (p^2 -3p + 6)/2$ occurrences of $a$ (the possible occurrences being at columns labeled by pairs of distinct entries containing neither $a$ nor $b$ and at a column labeled by $ab$). By symmetry,  also the number of occurrences of $b$ in such a row must be in the interval $[p - 1,(p^2 -3p + 6)/2]$. For an entry $c$ different from both $a$ and $b$ the number of occurrences in the row labeled by $ab$ is must be in the interval $[1,p-2]$ as the only forced occurrence of $c$ is in the column labeled by $cc$ and the possible occurrences in $p-3$ columns labeled by $cd$ for $d \not\in \{a,b,c\}$.

Again notice that intervals $[1,p-2]$, $[p - 1,(p^2 -3p + 6)/2]$, and $[(p^2 - p + 2)/2,p(p+1)/2]$ do not intersect (here we are using $p \geq 3$) and so we can use the entry count to determine the correct row labels. If the number of occurrences of the most frequent entry (say $a$) lies in $[(p^2 - p + 2)/2,p(p+1)/2]$, all other entries must occur at most once and the row label is $aa$. If the number of occurrences of the most frequent entry (say $a$) lies in $[p - 1,(p^2 -3p + 6)/2]$ then for some other entry (say $b$) the number of occurrences must lie in the same interval, for all other entries the count is in $[1,p-2]$, and the row label is $ab$. If any of the counts lies outside its prescribed interval or if the assigned labels do not exhaust all possible label (each exactly once), then the input matrix is not a DAV form. Otherwise, once we assign the row labels and check that all counts are in the right intervals, we proceed to assign the column labels in the same way as in Theorem~\ref{t2}.

\subsection{Case $p=2$ and arbitrary $n$.}
\label{p=2}
In this case we have $n\geq 2$ voters and each voter $i\in I = \{1, \ldots ,n\}$ has $\gal_i$ voting cards that can be distributed between the two candidates, denoted here by $A=\{a,b\}$. For a nonnegative integer $z$ the $2$-dimensional vector $(z,\gal_i-z)$ can be viewed as a strategy of voter $i\in I$ if $\gal_i-z\geq 0$. To simplify our notation, we simply use integers $z_i\in\ZZ_{\gal_i}=\{0,1,...,\gal_i\}$ to denote the corresponding unique strategy $(z_i,\gal_i-z_i)$ of voter $i\in I$. Thus, voting correspondences and forms must have dimensions $(\gal_1+1)\times \cdots \times (\gal_n+1)$ and $Z=\ZZ_{\gal_1}\times\dots\times\ZZ_{\gal_n}$ is the set of strategies of the voters. 

Let us first observe the following easy property:

\begin{lemma}\label{l-2n-1}
If $h:Z\to 2^A$ is a DAV correspondence, $g:Z\to A$ is a DAV form, and $\sigma=\sum_{i\in I}\gal_i$, then for all $(z_1,\dots,z_n)\in Z$ we have
\begin{equation}\label{e-2n-1}
h(z_1,\dots,z_n) ~=~ \begin{cases}
\{a\}&\text{for }\displaystyle\sum_{i\in I} z_i ~>~ \frac{\sigma}{2}\\
\{a,b\}&\text{for }\displaystyle\sum_{i\in I} z_i ~=~ \frac{\sigma}{2}\\
\{b\}&\text{for }\displaystyle\sum_{i\in I} z_i ~<~ \frac{\sigma}{2}
\end{cases}
\end{equation}
and
\begin{equation}\label{e-2n-2}
g(z_1,\dots,z_n) ~=~\begin{cases}
 a&\text{for }\displaystyle\sum_{i\in I} z_i ~>~ \frac{\sigma}{2}\\
a \text{ or } b&\text{for }\displaystyle\sum_{i\in I} z_i ~=~ \frac{\sigma}{2}\\
b&\text{for }\displaystyle\sum_{i\in I} z_i ~<~ \frac{\sigma}{2}
\end{cases}
\end{equation}
\qed
\end{lemma}
Note that if $\sigma$ is odd, then we have no difference between DAV correspondences and forms.

Given an $n$-dimensional array $f$, let  us denote by $f_{i,z}$ its $(n-1)$-dimensional subarray obtained by fixing the $i$th coordinate at value $z$, and varying all other coordinates over their possible values. We call $f_{i,z}$ a plane of $f$ (called a row or a column in the $2$-dimensional case). 

For a voter $i\in I$, integer $z\in \ZZ_{\gal_i}$, candidate $c\in A$, and voting correspondence $h: Z\to 2^A$ let us denote by
\[
s_h(i,z,c) ~=~ |\{(z_1,\dots,z_n)\in Z\mid z_i=z \text{ and } c\in h(z_1,\dots,z_n)\}|
\]
the number of occurrences of $c$ in the plane $h_{i,z}$. Similarly for a voting form $g:Z\to A$ let us denote by
\[
s_g(i,z,c) ~=~ |\{(z_1,\dots,z_n)\in Z\mid z_i=z \text{ and } c= g(z_1,\dots,z_n)\}|
\]
the number of occurrences of $c$ in the plane $g_{i,z}$. 

\begin{lemma}\label{l-2n-2}
If $h:Z\to 2^A$ is a DAV correspondence, $i\in I$, and $z,z'\in\ZZ_{\gal_i}$, $z\neq z'$ then we have
$(s_h(i,z,a),s_h(i,z,b))=(s_h(i,z',a),s_h(i,z',b))$ if and only if the planes $h_{i,z}$ and $h_{i,z'}$ are identical. Similarly, 
for a DAV form $g:Z\to A$ we have
$s_g(i,z,a)=s_g(i,z',a)$ if and only if the planes $g_{i,z}$ and $g_{i,z'}$ are identical. 
\end{lemma}

\proof
Immediate by Lemma \ref{l-2n-1}.
\qed

The above claims readily imply that DAV correspondences and forms can be recognized in polynomial time:

\begin{theorem}
Given a voting correspondence $h:Z\to 2^A$ let us relabel its planes such that for all $i\in I$ the sequence $(s_h(i,z,a)\mid z\in \ZZ_{\gal_i})$ is monotone non-increasing. Then $h$ is a DAV correspondence if and only if equalities \eqref{e-2n-1} hold for all $(z_1,\dots,z_n)\in Z$. Computing $s_h(i,z,a)$ for $i\in I$, $z\in \ZZ_{\gal_i}$, relabeling planes, and checking \eqref{e-2n-1} can be done in $O(\prod_{i\in I} (\gal_i+1))$ time.
\end{theorem}

\begin{theorem}
Given a voting form $g:Z\to A$ let us relabel its planes such that for all $i\in I$ the sequence $(s_g(i,z,a)\mid z\in \ZZ_{\gal_i})$ is monotone non-increasing. Then $g$ is a DAV form if and only if relations \eqref{e-2n-2} hold for all $(z_1,\dots,z_n)\in Z$. Computing $s_g(i,z,a)$ for $i\in I$, $z\in \ZZ_{\gal_i}$,  relabeling planes, and checking \eqref{e-2n-2} can be done in $O(\prod_{i\in I} (\gal_i+1))$ time.
\end{theorem}

\bigskip

\section{Conclusions}

In this paper we deal with distributed approval voting (DAV) correspondences and DAV forms for two voters ($n=2$). For both correspondences and forms we have two main results. The first result gives necessary and sufficient conditions under which every DAV correspodence (form) has distinct rows (and symmetrical conditions for distinct columns). The second result is a polynomial time recognition algorithm that for a given input voting correspondence (form) decides whether it is a DAV correspondence (form). In the affirmative case the algorithm also outputs a corresponding labelling of rows and columns by voting strategies. The main difference between correspondences and forms rests in the fact that while the recognition algorithm for DAV correspondences admits any input, the recognition algorithm for DAV forms has additional conditions on the input parameters, and thus does not give an answer for all possible inputs. The last section of the paper treats three special cases. The first special case - two voters, each having only one vote - is the original plurality voting (for two voters) which motivated the concept of distributed approval voting. For this case we give a complete characterization of DAV voting forms by a set of three forbidden submatrices. For the second special case - two voters, each having at most two votes - we give a simple polynomial time recognition algorithm for DAV forms as this case is not fully covered by the general recognition algorithm mentioned above. The last special case characterizes the case of two candidates and an arbitrary number of voters. 

There are two directions in which the problems solved in this paper could be generalized. One is to consider partial DAV correspondences and partial DAV forms. By "partial" we mean the following: both rows and columns are still labeled by voting strategies (and a voting strategy is still a particular distribution of of voting cards of a given voter among all candidates), but the set of labels may now be any subset of the set of all voting strategies. Recognizing such "subcorrespondences" and "subforms" in polynomial time seems to be hard to achieve. Another  generalization is to have simultaneously more than two candidates ($p > 2$) and more than two voters ($n > 2$). This second generalization of course works for both the original problems and for the partial ones. 

\section*{Acknowledgements}

The second author gratefully acknowledges a support by the Czech Science Foundation (Grant 19-19463S),
the third author was partially supported by the RSF grant  20-11-20203.


\begin{thebibliography}{99}

\bibitem{Alo06}
C. Alos-Ferrer,  A simple characterization of approval voting,
Social Choice and Welfare 27 (2006) 621--625

\bibitem{Berge}
C. Berge, Two theorems in graph theory,
Proceedings of the National Academy of Sciences of the United States of America, 43 (9)  (1957)  842–844.

\bibitem{BCG19}
E. Boros, O. Cepek, and V. Gurvich.
Separable discrete functions: recognition and sufficient conditions;
Discrete Mathematics 342 (2019) 1275--1292.


\bibitem{BCG11}
E.~Boros, O.~\v{C}epek, and V.~Gurvich,
Total tightness implies Nash-solvability
for three-person game forms,
Discrete Mathematics 312:8 (2012) 1436--1443.




\bibitem{BGMP10}
E. Boros, V. Gurvich, K. Makino, and D. Papp,
Acyclic, or totally tight, two-person game forms;
a characterization and main properties;
Discrete Mathematics, 310 (6-7) (2010) 1135--1151.









\bibitem{BF78}
S. Brams and P. Fishburn,  Approval voting,
American Political Science Review. 72:3 (1978) 831--847.

\bibitem{BF07}
S. Brams and P. Fishburn,  Approval voting, Springer-Verlag, 2007.

\bibitem{DL04}
A. Dhillon and B. Lockwood,
When are pluralityrule voting games dominance-solvable?
Games and Economic Behavior 46 (2004) 55--75.

\bibitem{Dau02}
I. Daubechies, 
Weighted Voting Systems, Voting and Social Choice, Math Alive,
Princeton University, 2002. 

%
%






\bibitem{Egervary31}
J. Egerv{\'a}ry, Matrixok kombinatorius tulajdons{\'a}gair{\'o}l [Hungarian, with German summary], \emph{Matematikai {\'e}s Fizikai Lapok} 38 (1931) 16–28 [English translation [by H.W.
Kuhn]: On combinatorial properties of matrices, Logistics Papers, George Washington University, issue 11 (1955), paper 4, pp. 1–11].


\bibitem{Konig31}
D. K{\"o}nig, Graphok {\'e}s matrixok [Hungarian; Graphs and matrices], \emph{Matematikai {\'e}s
Fizikai Lapok} 38 (1931) 116–119.

\bibitem{Kuhn55}
H.W. Kuhn, The Hungarian method for the assignment problem, \emph{Naval Research
Logistics Quarterly} 2 (1955) 83–97.









\bibitem{Kuk11}
N.S. Kukushkin
Acyclicity of improvements in finite game forms,
International Journal of Game Theory 40 (2011) 147--177.

\bibitem{MV06}
J. Masso and M. Vorsatz, 
Weighted Approval Voting,
Economic Theory 36 (2008) 129--146.

\bibitem{Mei15}
R. Meir, Plurality voting under uncertainty,
Proc. of 29th AAAI (2015) 2103--2109.

\bibitem{MPRJ17}
R. Meir, M. Polukarov, J.S. Rosenschein, and N.R. Jennings,
Iterative voting and acyclic games,
Artificial Intelligence 252 (2017) 100--122.

\bibitem{Mou79}
H.Moulin, Dominance-solvable voting schemes, Econometrica 47 (1979) 1337-1351.

\bibitem{Mou83}
H.Moulin, The strategy of social choice.
North-Holland Publ. Co, Amsterdam, New York, Oxford, 1983.

\bibitem{Pel84}
B. Peleg, Game theoretic analysis of voting in committees,
Cambridge University Press, 1984.


 \bibitem{Sen70}
 A. K. Sen, Collective Choice and Social Welfare, Holden-Day, 1970.

\end{thebibliography}
\end{document}